\documentclass[12pt]{article}
\usepackage{amsmath,amssymb,amsthm,latexsym,amscd}

\title{Approximations of theories\footnote{
\newline\indent \ \ \ This research was partially
supported by Russian Foundation for Basic Researches (Project No.
17-01-00531-a) and Committee of Science in Education and Science
Ministry of the Republic of Kazakhstan (Grant No. AP05132546).} }
\author{Sergey V.
Sudoplatov\footnote{sudoplat@math.nsc.ru}}
\date{}
\begin{document}
\maketitle

\begin{abstract}
We study approximations of theories both in general context and
with respect to some natural classes of theories. Some kinds of
approximations are considered, connections with finitely
axiomatizable theories and minimal generating sets of theories as
well as their $e$-spectra are found.

\end{abstract}

Approximations of structures and theories as well as closures with
respect to these approximations were studied in a series of
papers, both implicitly \cite{Bank1, Bank2, ChHL, ChL, AZ, Rosen,
Vaa, ChHr, MS, AK} and explicitly  \cite{cs, cl, lut, ccct, rest,
lft}. They are connected with the technique for finitely
axiomatizable theories \cite{EfrFuh, Pa1, Mak1, Pa2, Zilb, Hr94,
Peret, BlBou}.

The aim of the paper is to introduce and investigate
approximations of theories both in general context and with
respect to some natural classes of theories.

The paper is organized as follows. In Section 1 we collect
preliminary definitions and assertions. In Section 2 we define
approximations relative given family $\mathcal{T}$ of theories and
characterize the property ``to be $\mathcal{T}$-approximated''. In
Section 3 we connect approximable theories with finite
axiomatizable ones, introduce the notion of
$\mathcal{T}$-relatively finite axiomatizability and characterize
this notion. In Section 4 we consider $\lambda$-approximable
theories, i.e., theories generated by families of theories such
that these families have the cardinality $\lambda$, and
characterize the property of $\lambda$-approximation.
Approximations by almost language uniform theories are considered
in Section 5. A characterization for approximating subfamilies and
lower bounds for $e$-spectra are proved in Section 6. In Section
7, $e$-categorical approximating families are introduced and
characterized.

\section{Preliminaries}

Throughout the paper we consider complete first-order theories $T$
in predicate languages $\Sigma(T)$ and use the following
terminology in \cite{cs, cl, lut, ccct, rest, lft}.

Let $P=(P_i)_{i\in I}$, be a family of nonempty unary predicates,
$(\mathcal{A}_i)_{i\in I}$ be a family of structures such that
$P_i$ is the universe of $\mathcal{A}_i$, $i\in I$, and the
symbols $P_i$ are disjoint with languages for the structures
$\mathcal{A}_j$, $j\in I$. The structure
$\mathcal{A}_P\rightleftharpoons\bigcup\limits_{i\in
I}\mathcal{A}_i$\index{$\mathcal{A}_P$} expanded by the predicates
$P_i$ is the {\em $P$-union}\index{$P$-union} of the structures
$\mathcal{A}_i$, and the operator mapping $(\mathcal{A}_i)_{i\in
I}$ to $\mathcal{A}_P$ is the {\em
$P$-operator}\index{$P$-operator}. The structure $\mathcal{A}_P$
is called the {\em $P$-combination}\index{$P$-combination} of the
structures $\mathcal{A}_i$ and denoted by ${\rm
Comb}_P(\mathcal{A}_i)_{i\in I}$\index{${\rm
Comb}_P(\mathcal{A}_i)_{i\in I}$} if
$\mathcal{A}_i=(\mathcal{A}_P\upharpoonright
A_i)\upharpoonright\Sigma(\mathcal{A}_i)$, $i\in I$. Structures
$\mathcal{A}'$, which are elementary equivalent to ${\rm
Comb}_P(\mathcal{A}_i)_{i\in I}$, will be also considered as
$P$-combinations.

Clearly, all structures $\mathcal{A}'\equiv {\rm
Comb}_P(\mathcal{A}_i)_{i\in I}$ are represented as unions of
their restrictions $\mathcal{A}'_i=(\mathcal{A}'\upharpoonright
P_i)\upharpoonright\Sigma(\mathcal{A}_i)$ if and only if the set
$p_\infty(x)=\{\neg P_i(x)\mid i\in I\}$ is inconsistent. If
$\mathcal{A}'\ne{\rm Comb}_P(\mathcal{A}'_i)_{i\in I}$, we write
$\mathcal{A}'={\rm Comb}_P(\mathcal{A}'_i)_{i\in
I\cup\{\infty\}}$, where
$\mathcal{A}'_\infty=\mathcal{A}'\upharpoonright
\bigcap\limits_{i\in I}\overline{P_i}$, maybe applying
Morleyzation. Moreover, we write ${\rm
Comb}_P(\mathcal{A}_i)_{i\in I\cup\{\infty\}}$\index{${\rm
Comb}_P(\mathcal{A}_i)_{i\in I\cup\{\infty\}}$} for ${\rm
Comb}_P(\mathcal{A}_i)_{i\in I}$ with the empty structure
$\mathcal{A}_\infty$.

Note that if all predicates $P_i$ are disjoint, a structure
$\mathcal{A}_P$ is a $P$-combination and a disjoint union of
structures $\mathcal{A}_i$. In this case the $P$-combination
$\mathcal{A}_P$ is called {\em
disjoint}.\index{$P$-combination!disjoint} Clearly, for any
disjoint $P$-combination $\mathcal{A}_P$, ${\rm
Th}(\mathcal{A}_P)={\rm Th}(\mathcal{A}'_P)$, where
$\mathcal{A}'_P$ is obtained from $\mathcal{A}_P$ replacing
$\mathcal{A}_i$ by pairwise disjoint
$\mathcal{A}'_i\equiv\mathcal{A}_i$, $i\in I$. Thus, in this case,
similar to structures the $P$-operator works for the theories
$T_i={\rm Th}(\mathcal{A}_i)$ producing the theory $T_P={\rm
Th}(\mathcal{A}_P)$\index{$T_P$}, being {\em
$P$-combination}\index{$P$-combination} of $T_i$, which is denoted
by ${\rm Comb}_P(T_i)_{i\in I}$.\index{${\rm Comb}_P(T_i)_{i\in
I}$}

For an equivalence relation $E$ replacing disjoint predicates
$P_i$ by $E$-classes we get the structure
$\mathcal{A}_E$\index{$\mathcal{A}_E$} being the {\em
$E$-union}\index{$E$-union} of the structures $\mathcal{A}_i$. In
this case the operator mapping $(\mathcal{A}_i)_{i\in I}$ to
$\mathcal{A}_E$ is the {\em $E$-operator}\index{$E$-operator}. The
structure $\mathcal{A}_E$ is also called the {\em
$E$-combination}\index{$E$-combination} of the structures
$\mathcal{A}_i$ and denoted by ${\rm Comb}_E(\mathcal{A}_i)_{i\in
I}$\index{${\rm Comb}_E(\mathcal{A}_i)_{i\in I}$}; here
$\mathcal{A}_i=(\mathcal{A}_E\upharpoonright
A_i)\upharpoonright\Sigma(\mathcal{A}_i)$, $i\in I$. Similar
above, structures $\mathcal{A}'$, which are elementary equivalent
to $\mathcal{A}_E$, are denoted by ${\rm
Comb}_E(\mathcal{A}'_j)_{j\in J}$, where $\mathcal{A}'_j$ are
restrictions of $\mathcal{A}'$ to its $E$-classes. The
$E$-operator works for the theories $T_i={\rm Th}(\mathcal{A}_i)$
producing the theory $T_E={\rm Th}(\mathcal{A}_E)$\index{$T_E$},
being {\em $E$-combination}\index{$E$-combination} of $T_i$, which
is denoted by ${\rm Comb}_E(T_i)_{i\in I}$\index{${\rm
Comb}_E(T_i)_{i\in I}$} or by ${\rm
Comb}_E(\mathcal{T})$\index{${\rm Comb}_E(\mathcal{T})$}, where
$\mathcal{T}=\{T_i\mid i\in I\}$.

Clearly, $\mathcal{A}'\equiv\mathcal{A}_P$ realizing $p_\infty(x)$
is not elementary embeddable into $\mathcal{A}_P$ and can not be
represented as a disjoint $P$-combination of
$\mathcal{A}'_i\equiv\mathcal{A}_i$, $i\in I$. At the same time,
there are $E$-combinations such that all
$\mathcal{A}'\equiv\mathcal{A}_E$ can be represented as
$E$-combinations of some $\mathcal{A}'_j\equiv\mathcal{A}_i$. We
call this representability of $\mathcal{A}'$ to be the {\em
$E$-representability}.

If there is $\mathcal{A}'\equiv\mathcal{A}_E$ which is not
$E$-representable, we have the $E'$-representability replacing $E$
by $E'$ such that $E'$ is obtained from $E$ adding equivalence
classes with models for all theories $T$, where $T$ is a theory of
a restriction $\mathcal{B}$ of a structure
$\mathcal{A}'\equiv\mathcal{A}_E$ to some $E$-class and
$\mathcal{B}$ is not elementary equivalent to the structures
$\mathcal{A}_i$. The resulting structure $\mathcal{A}_{E'}$ (with
the $E'$-representability) is a {\em
$e$-completion}\index{$e$-completion}, or a {\em
$e$-saturation}\index{$e$-saturation}, of $\mathcal{A}_{E}$. The
structure $\mathcal{A}_{E'}$ itself is called {\em
$e$-complete}\index{Structure!$e$-complete}, or {\em
$e$-saturated}\index{Structure!$e$-saturated}, or {\em
$e$-universal}\index{Structure!$e$-universal}, or {\em
$e$-largest}\index{Structure!$e$-largest}.

For a structure $\mathcal{A}_E$ the number of {\em
new}\index{Structure!new} structures with respect to the
structures $\mathcal{A}_i$, i.~e., of the structures $\mathcal{B}$
which are pairwise elementary non-equivalent and elementary
non-equivalent to the structures $\mathcal{A}_i$, is called the
{\em $e$-spectrum}\index{$e$-spectrum} of $\mathcal{A}_E$ and
denoted by $e$-${\rm Sp}(\mathcal{A}_E)$.\index{$e$-${\rm
Sp}(\mathcal{A}_E)$} The value ${\rm sup}\{e$-${\rm
Sp}(\mathcal{A}'))\mid\mathcal{A}'\equiv\mathcal{A}_E\}$ is called
the {\em $e$-spectrum}\index{$e$-spectrum} of the theory ${\rm
Th}(\mathcal{A}_E)$ and denoted by $e$-${\rm Sp}({\rm
Th}(\mathcal{A}_E))$.\index{$e$-${\rm Sp}({\rm
Th}(\mathcal{A}_E))$} If structures $\mathcal{A}_i$ represent
theories $T_i$ of a family $\mathcal{T}$, consisting of $T_i$,
$i\in I$, then the $e$-spectrum $e$-${\rm Sp}(\mathcal{A}_E)$ is
denoted by $e$-${\rm Sp}(\mathcal{T})$.

If $\mathcal{A}_E$ does not have $E$-classes $\mathcal{A}_i$,
which can be removed, with all $E$-classes
$\mathcal{A}_j\equiv\mathcal{A}_i$, preserving the theory ${\rm
Th}(\mathcal{A}_E)$, then $\mathcal{A}_E$ is called {\em
$e$-prime}\index{Structure!$e$-prime}, or {\em
$e$-minimal}\index{Structure!$e$-minimal}.

For a structure $\mathcal{A}'\equiv\mathcal{A}_E$ we denote by
${\rm TH}(\mathcal{A}')$ the set of all theories ${\rm
Th}(\mathcal{A}_i)$\index{${\rm Th}(\mathcal{A}_i)$} of
$E$-classes $\mathcal{A}_i$ in $\mathcal{A}'$.

By the definition, an $e$-minimal structure $\mathcal{A}'$
consists of $E$-classes with a minimal set ${\rm
TH}(\mathcal{A}')$. If ${\rm TH}(\mathcal{A}')$ is the least for
models of ${\rm Th}(\mathcal{A}')$ then $\mathcal{A}'$ is called
{\em $e$-least}.\index{Structure!$e$-least}

\medskip
{\bf Definition} \cite{cl}. Let $\overline{\mathcal{T}}_\Sigma$ be
the set of all complete elementary theories of a relational
language $\Sigma$. For a set
$\mathcal{T}\subset\overline{\mathcal{T}}_\Sigma$ we denote by
${\rm Cl}_E(\mathcal{T})$ the set of all theories ${\rm
Th}(\mathcal{A})$, where $\mathcal{A}$ is a structure of some
$E$-class in $\mathcal{A}'\equiv\mathcal{A}_E$,
$\mathcal{A}_E={\rm Comb}_E(\mathcal{A}_i)_{i\in I}$, ${\rm
Th}(\mathcal{A}_i)\in\mathcal{T}$. As usual, if $\mathcal{T}={\rm
Cl}_E(\mathcal{T})$ then $\mathcal{T}$ is said to be {\em
$E$-closed}.\index{Set!$E$-closed}

The operator ${\rm Cl}_E$ of $E$-closure can be naturally extended
to the classes $\mathcal{T}\subset\overline{\mathcal{T}}$, where
$\overline{\mathcal{T}}$ is the union of all
$\overline{\mathcal{T}}_\Sigma$ as follows: ${\rm
Cl}_E(\mathcal{T})$ is the union of all ${\rm
Cl}_E(\mathcal{T}_0)$ for subsets
$\mathcal{T}_0\subseteq\mathcal{T}$, where new language symbols
with respect to the theories in $\mathcal{T}_0$ are empty.

For a set $\mathcal{T}\subset\overline{\mathcal{T}}$ of theories
in a language $\Sigma$ and for a sentence $\varphi$ with
$\Sigma(\varphi)\subseteq\Sigma$ we denote by
$\mathcal{T}_\varphi$\index{$\mathcal{T}_\varphi$} the set
$\{T\in\mathcal{T}\mid\varphi\in T\}$. Any set
$\mathcal{T}_\varphi$ is called a {\em $\varphi$-neighbourhood},
or simply a {\em neighbourhood}, for $\mathcal{T}$.

\medskip
{\bf Proposition 1.1} \cite{cl}. {\em If
$\mathcal{T}\subset\overline{\mathcal{T}}$ is an infinite set and
$T\in\overline{\mathcal{T}}\setminus\mathcal{T}$ then $T\in{\rm
Cl}_E(\mathcal{T})$ {\rm (}i.e., $T$ is an {\sl accumulation
point} for $\mathcal{T}$ with respect to $E$-closure ${\rm
Cl}_E${\rm )} if and only if for any formula $\varphi\in T$ the
set $\mathcal{T}_\varphi$ is infinite.}

\medskip
If $T$ is an accumulation point for $\mathcal{T}$ then we also say
that $T$ is an {\em accumulation point} for ${\rm
Cl}_E(\mathcal{T})$.

\medskip
{\bf Theorem 1.2} \cite{cl}. {\em For any sets
$\mathcal{T}_0,\mathcal{T}_1\subset\overline{\mathcal{T}}$, ${\rm
Cl}_E(\mathcal{T}_0\cup\mathcal{T}_1)={\rm
Cl}_E(\mathcal{T}_0)\cup{\rm Cl}_E(\mathcal{T}_1)$.}

\medskip
{\bf Definition} \cite{cl}. Let $\mathcal{T}_0$ be a closed set in
a topological space $(\mathcal{T},\mathcal{O}_E(\mathcal{T}))$,
where $\mathcal{O}_E(\mathcal{T})=\{\mathcal{T}\setminus{\rm
Cl}_E(\mathcal{T}')\mid\mathcal{T}'\subseteq\mathcal{T}\}$. A
subset $\mathcal{T}'_0\subseteq\mathcal{T}_0$ is said to be {\em
generating}\index{Set!generating} if $\mathcal{T}_0={\rm
Cl}_E(\mathcal{T}'_0)$. The generating set $\mathcal{T}'_0$ (for
$\mathcal{T}_0$) is {\em minimal}\index{Set!generating!minimal} if
$\mathcal{T}'_0$ does not contain proper generating subsets. A
minimal generating set $\mathcal{T}'_0$ is {\em
least}\index{Set!generating!least} if $\mathcal{T}'_0$ is
contained in each generating set for $\mathcal{T}_0$.

\medskip
{\bf Theorem 1.3} \cite{cl}. {\em If $\mathcal{T}'_0$ is a
generating set for a $E$-closed set $\mathcal{T}_0$ then the
following conditions are equivalent:

$(1)$ $\mathcal{T}'_0$ is the least generating set for
$\mathcal{T}_0$;

$(2)$ $\mathcal{T}'_0$ is a minimal generating set for
$\mathcal{T}_0$;

$(3)$ any theory in $\mathcal{T}'_0$ is isolated by some set
$(\mathcal{T}'_0)_\varphi$, i.e., for any $T\in\mathcal{T}'_0$
there is $\varphi\in T$ such that
$(\mathcal{T}'_0)_\varphi=\{T\}$;

$(4)$ any theory in $\mathcal{T}'_0$ is isolated by some set
$(\mathcal{T}_0)_\varphi$, i.e., for any $T\in\mathcal{T}'_0$
there is $\varphi\in T$ such that
$(\mathcal{T}_0)_\varphi=\{T\}$.}

\medskip
{\bf Proposition 1.4.} {\em Any family $\mathcal{T}$ of theories
can be expanded till a family $\mathcal{T}'$ with the least
generating set.}

\medskip
Proof. By Theorem 1.3 it suffices to introduce, for each theory
$T\in \mathcal{T}$, a new unary predicate $P_T$ such that $P_T$ is
complete for $T$ and empty for all $T'\in
\mathcal{T}\setminus\{T\}$. Clearly, that the formula witnessing
that $P_T$ is complete separates $T$ from
$\mathcal{T}\setminus\{T\}$. Thus, the family $\mathcal{T}'$
itself is the least generating set.~$\Box$

\section{$\mathcal{T}$-approximations}

{\bf Definition.} Let $\mathcal{T}$ be a class of theories and $T$
be a theory, $T\notin\mathcal{T}$. The theory $T$ is called {\em
$\mathcal{T}$-approximated}, or {\em approximated by}
$\mathcal{T}$, or {\em $\mathcal{T}$-approximable}, or a {\em
pseudo-$\mathcal{T}$-theory}, if for any formula $\varphi\in T$
there is $T'\in\mathcal{T}$ such that $\varphi\in T'$.

If $T$ is $\mathcal{T}$-approximated then $\mathcal{T}$ is called
an {\em approximating family} for $T$, and theories
$T'\in\mathcal{T}$ are {\em approximations} for $T$.

\medskip
{\bf Remark 2.1.} If $T$ is $\mathcal{T}$-approximated then $T$ is
$\mathcal{T}'$-approximated for any
$\mathcal{T}'\supseteq\mathcal{T}$. At the same time, if $T$ is
$\mathcal{T}$-approximated then $T$ is
$\mathcal{T}\setminus\{T'\}$-approximated for any
$T'\in\mathcal{T}$. Indeed, since $T'\ne T$, there is $\psi\in T$
such that $\neg\psi\in T'$, and for any formula $\varphi\in T$ the
formula $\varphi\wedge\psi$ belongs both to $T$ and to some
$T''\in \mathcal{T}\setminus\{T'\}$, so $\varphi\in T''$.

Besides, an approximation family $\mathcal{T}$ for $T$ can be
extended by an arbitrary theory $T'\ne T$, assuming the
possibility to extend the language $\Sigma(T)$. Thus, if there an
approximating family $\mathcal{T}$ for $T$ then $\mathcal{T}$ can
not be chosen minimal or maximal by inclusion, and if the language
$\Sigma(T)$ is fixed then the maximal one exists containing all
$\Sigma(T)$-theories $T'\ne T$.

\medskip
Remark 2.1 implies the following proposition, but we will give
slightly different arguments.

\medskip
{\bf Proposition 2.1.} {\em If there is a
$\mathcal{T}$-approximated theory then $\mathcal{T}$ is infinite.}

\medskip
Proof. If $T$ is $\mathcal{T}$-approximated and $\mathcal{T}$ is
finite consisting of $T_1,\ldots,T_n$ then having
$T\notin\mathcal{T}$ there are formulas $\varphi_i\in T_i$ such
that
$\psi\rightleftharpoons\neg\varphi_1\wedge\ldots\wedge\neg\varphi_n\in
T$. Since $\psi\notin T_1\cup\ldots\cup T_n$, $T$ can not be
$\mathcal{T}$-approximated, implying a contradiction. $\Box$

\medskip
{\bf Proposition 2.2.} {\em A theory $T\notin\mathcal{T}$ is
$\mathcal{T}$-approximated if and only if $T\in{\rm
Cl}_E(\mathcal{T})$.}

\medskip
Proof. Let $T$ be $\mathcal{T}$-approximated. By Proposition 1.1
it suffices to show that for any $\varphi\in T$ there are
infinitely many $T'\in\mathcal{T}$ such that $\varphi\in T'$.
Assuming on contrary that there are only finitely many $T'$, say
$T_1,\ldots,T_n$, then there are $\varphi_i\in T_i$ such that
$\psi\rightleftharpoons\varphi\wedge\neg\varphi_1\wedge\ldots\wedge\neg\varphi_n\in
T$. Since $\psi$ does not belong to $\cup\mathcal{T}$, then $T$ is
not $\mathcal{T}$-approximated.

If $T\in{\rm Cl}_E(\mathcal{T})$ then, by Proposition 1.1, for any
$\varphi\in T$ there are infinitely many $T'\in\mathcal{T}$ such
that $\varphi\in T'$.~$\Box$

\medskip
Recall that $\mathcal{T}$ is {\em $E$-closed} if $\mathcal{T}={\rm
Cl}_E(\mathcal{T})$. By Proposition 2.2 we have

\medskip
{\bf Corollary 2.3.} {\em For any family $\mathcal{T}$ there is a
$\mathcal{T}$-approximated theory if and only if $\mathcal{T}$ is
not $E$-closed.}

\medskip
{\bf Definition} \cite{Rosen}. An infinite structure $\mathcal{M}$
is {\em pseudo-finite} if every sentence true in $\mathcal{M}$ has
a finite model.

\medskip
If $T={\rm Th}(\mathcal{M})$ for pseudo-finite $\mathcal{M}$ then
$T$ is called {\em pseudo-finite} as well.

Following \cite{ccct} we denote by $\overline{\mathcal{T}}$ the
class of all complete elementary theories of relational languages,
by $\overline{\mathcal{T}}_{\rm fin}$ the subclass of
$\overline{\mathcal{T}}$ consisting of all theories with finite
models, and by $\overline{\mathcal{T}}_{\rm inf}$ the class
$\overline{\mathcal{T}}\setminus\overline{\mathcal{T}}_{\rm fin}$.

\medskip
{\bf Proposition 2.4.} {\em For any theory $T$ the following
conditions are equivalent:

$(1)$ $T$ is pseudo-finite;

$(2)$ $T$ is $\overline{\mathcal{T}}_{\rm fin}$-approximated;

$(3)$ $T\in{\rm Cl}_E(\overline{\mathcal{T}}_{\rm
fin})\setminus\overline{\mathcal{T}}_{\rm fin}$.}

\medskip
Proof. $(1)\Leftrightarrow(2)$ holds by the definition.
$(2)\Leftrightarrow(3)$ is satisfied by Proposition 2.2.~$\Box$

\section{Approximable and finitely axiomatizable theories}

{\bf Definition.} A theory $T$ is called {\em approximable} if $T$
is $\mathcal{T}$-approximable for some $\mathcal{T}$.

Recall \cite{Peret} that a theory $T$ is {\em finitely
axiomatizable} if $T$ is forced by some formula $\varphi\in T$.

Notice that by the definition finitely axiomatizable theories have
finite languages.

\medskip
{\bf Proposition 3.1.} {\em For any theory $T$ the following
conditions are equivalent:

$(1)$ $T$ is approximable;

$(2)$ $T$ is $\overline{\mathcal{T}}\setminus\{T\}$-approximated;

$(3)$ $T$ is not finitely axiomatizable.}

\medskip
Proof. $(1)\Leftrightarrow(2)$ holds by the definition.

$(2)\Rightarrow(3)$. Assume that $T$ is finitely axiomatizable
witnessed by a formula $\varphi\in T$. Then $\varphi\notin T'$ for
any $T'\in\overline{\mathcal{T}}\setminus\{T\}$. Hence, $T$ is not
$\mathcal{T}$-approximated for any $\mathcal{T}\subseteq
\overline{\mathcal{T}}$. In particular, $T$ in not
$\overline{\mathcal{T}}\setminus\{T\}$-approximated.

$(3)\Rightarrow(2)$. Let $T$ be not finitely axiomatizable. Then
for any $\varphi\in T$ there is
$T'\in\overline{\mathcal{T}}\setminus\{T\}$ with $\varphi\in T'$,
since otherwise $T$ is axiomatizable by $\varphi$. Therefore, $T$
is $\overline{\mathcal{T}}\setminus\{T\}$-approximated.~$\Box$

\medskip
Illustrating Proposition 3.1 we note that any theory $T$, in an
infinite relational language $\Sigma$, is approximable by theories
$T\upharpoonright\Sigma_0$, for finite $\Sigma_0$, expanded by
empty or complete predicates for symbols
$P\in\Sigma\setminus\Sigma_0$ such that $P$ is empty for these
expansions if $P$ is not empty for $T$, and $P$ is complete for
these expansions if $P$ is empty for $T$.

\medskip
We denote by $\overline{\mathcal{T}}_{\rm fa}$ the class of all
finitely axiomatizable theories, which coincide, by Proposition
3.1, with the class of all non-approximable theories. By the
definition the class $\overline{\mathcal{T}}_{\rm fa}$ consists
exactly of theories $T$ having singletons
$\overline{\mathcal{T}}_\varphi=\{T\}$, where $\varphi\vdash T$.
Thus, by Propositions 1.1, the class
$\overline{\mathcal{T}}\setminus\overline{\mathcal{T}}_{\rm fa}$
is $E$-closed, whereas $\overline{\mathcal{T}}_{\rm fa}$ is not
$E$-closed, whose $E$-closure contains at least all pseudo-finite
theories.

In the connection with this property it is natural to pose the
following

\medskip
{\bf Problem 1.} {\em Describe $\mathcal{T}$-approximable theories
and ${\rm Cl}_E(\mathcal{T})$ for natural classes
$\mathcal{T}\subseteq\overline{\mathcal{T}}_{\rm fa}$.}

\medskip
This problem can be considered in the following context.

\medskip
{\bf Definition.} For a family $\mathcal{T}$, a theory $T$ is {\em
$\mathcal{T}$-finitely axiomatizable}, or {\em finitely
axiomatizable with respect to $\mathcal{T}$}, or {\em
$\mathcal{T}$-relatively finitely axiomatizable}, if
$\mathcal{T}_\varphi=\{T\}$ for some
$\Sigma(\mathcal{T})$-sentence $\varphi$.

\medskip
{\bf Remark 3.2.} 1. A theory $T$ is finitely axiomatizable if and
only if $T$ is $\mathcal{T}$-finitely axiomatizable for any
$\mathcal{T}$ in the language $\Sigma(T)$.

2. A theory $T$ is $\mathcal{T}$-finitely axiomatizable if and
only if there is finite $\mathcal{T}_\varphi$ containing $T$, for
some $\Sigma(\mathcal{T})$-sentence $\varphi$.

\medskip
In this context Theorem 1.3 can be reformulated in the following
way.

\medskip
{\bf Theorem 3.3.} {\em If $\mathcal{T}'_0$ is a generating set
for a $E$-closed set $\mathcal{T}_0$ then the following conditions
are equivalent:

$(1)$ $\mathcal{T}'_0$ is the least generating set for
$\mathcal{T}_0$;

$(2)$ $\mathcal{T}'_0$ is a minimal generating set for
$\mathcal{T}_0$;

$(3)$ any theory in $\mathcal{T}'_0$ is $\mathcal{T}'_0$-finitely
axiomatizable;

$(4)$ any theory in $\mathcal{T}'_0$ is $\mathcal{T}_0$-finitely
axiomatizable.}

\medskip
Problem 1 can be divide into two possibilities with respect to
${\rm Cl}_E(\mathcal{T})$: with or without the least generating
sets. In particular, it admits the following refinement.

\medskip
{\bf Problem ${\bf 1'}$.} {\em Describe $\mathcal{T}$-approximable
theories and ${\rm Cl}_E(\mathcal{T})$ for natural sets
$\mathcal{T}$ containing $\mathcal{T}$-finitely axiomatizable
generating sets.}

\medskip
{\bf Definition.} For a family $\mathcal{T}$ of a language
$\Sigma$, a sentence $\varphi$ of the language $\Sigma$ is called
{\em $\mathcal{T}$-complete} if $\varphi$ isolates a unique theory
in $\mathcal{T}$, i.e., $\mathcal{T}_\varphi$ is a singleton.

\medskip
Clearly, a sentence $\varphi$ is complete if and only if $\varphi$
is $\mathcal{T}$-complete for any family $\mathcal{T}$ with a
theory $T\in\mathcal{T}$ containing $\varphi$.

Obviously, if $|\mathcal{T}_\varphi|\in\omega\setminus\{0\}$ then
each theory in $\mathcal{T}_\varphi$ contains a
$\mathcal{T}$-complete sentence, but not vice versa. Indeed,
$\mathcal{T}$ can consist of infinitely many finitely
axiomatizable theories, so each theory in $\mathcal{T}_{\forall
x(x\approx x)}$ contains a $\mathcal{T}$-complete sentence whereas
$|\mathcal{T}_{\forall x(x\approx x)}|\geq\omega$.

Since $\mathcal{T}$-complete sentences confirm the
$\mathcal{T}$-finite axiomatizability of theories in
$\mathcal{T}$, and a theory $T$ contains a $\mathcal{T}$-complete
sentence if and only if $T$ contains a disjunction of
$\mathcal{T}$-complete sentences, Theorem 3.3 admits the following
reformulation, with a slight extension.

\medskip
{\bf Theorem 3.4.} {\em If $\mathcal{T}'_0$ is a generating set
for a $E$-closed set $\mathcal{T}_0$ then the following conditions
are equivalent:

$(1)$ $\mathcal{T}'_0$ is the least generating set for
$\mathcal{T}_0$;

$(2)$ $\mathcal{T}'_0$ is a minimal generating set for
$\mathcal{T}_0$;

$(3)$ any theory in $\mathcal{T}'_0$ contains a
$\mathcal{T}'_0$-complete sentence;

$(4)$ any theory in $\mathcal{T}'_0$ contains a
$\mathcal{T}_0$-complete sentence;

$(5)$ any theory in $\mathcal{T}'_0$ contains a disjunction of
$\mathcal{T}'_0$-complete sentences;

$(6)$ any theory in $\mathcal{T}'_0$ contains a disjunction of
$\mathcal{T}_0$-complete sentences.}

\section{$\lambda$-approximable theories}

Below we consider a series of notions related to cardinalities for
approximations of theories.

\medskip
{\bf Definition.} Let $\lambda$ be a cardinality, $\mathcal{T}$ be
a family of theories. A theory $T$ is called {\em
$(\lambda,\mathcal{T})$-approximable}, or {\em
$\lambda$-approximable} ({\em $\lambda$-approximated}) by
$\mathcal{T}$, if $T$ is $\mathcal{T}'$-approximable for some
$\mathcal{T}'\subseteq\mathcal{T}$ with $|\mathcal{T}'|=\lambda$.
A theory $T$ is called {\em somewhere {\rm (}almost everywhere{\rm
)} $(\lambda,\mathcal{T})$-approximable} if
$T\upharpoonright\Sigma$ is
$(\lambda,\mathcal{T}\upharpoonright\Sigma)$-approximable for some
(any) $\Sigma\subseteq\Sigma(T)$, where $|\Sigma|+\omega=\lambda$
and $\mathcal{T}\upharpoonright\Sigma$ is the restriction of
theories in $\mathcal{T}$ till the language $\Sigma$.

A theory $T$ is called {\em exactly} (somewhere / almost
everywhere) $(\lambda,\mathcal{T})$-approxi\-m\-able if $T$ is
(somewhere / almost everywhere)
$(\lambda,\mathcal{T})$-approximable and $T$ is not (somewhere /
almost everywhere) $(\mu,\mathcal{T})$-approximable for
$\mu<\lambda$.

A theory $T$ is called {\em {\rm (}exactly / somewhere / almost
everywhere{\rm )} $\lambda$-approximable} if $T$ is (exactly /
somewhere / almost everywhere)
$(\lambda,\mathcal{T})$-approximable for some $\mathcal{T}$.

\medskip
{\bf Remark 4.1.} By the definition, if $T$ is (exactly /
somewhere) $\lambda$-approximable then $\lambda\geq\omega$.
Besides, if $T$ is almost everywhere
$(\omega,\mathcal{T})$-approximable then $\mathcal{T}$ has
infinitely many theories of structures having distinct finite
cardinalities, since restrictions to the empty language can be
approximated only by theories of finite structures. Moreover, by
Proposition 3.1, $T$ does not have finitely axiomatizable
restrictions.

Again by the definition, if $T$ is almost everywhere
$(\lambda,\mathcal{T})$-approximable then $T$ is almost
$(\lambda,\mathcal{T})$-approximable. But not vice versa, since
$T$ can contain finitely axiomatizable restrictions.

\medskip
If $\lambda=\omega$ we also say about {\em countably approximable
theories} instead of $\omega$-approximable and
$(\omega,\mathcal{T})$-approximable ones.

\medskip
Clearly by the definition that countably approximable theories are
exactly countably approximable.

The following problem is related to the series of notions above.

\medskip
{\bf Problem 2.} {\em Describe cardinalities $\lambda$ and forms
of approximations for natural classes of theories.}

\medskip
Illustrating the notions above and partially answering the problem
we consider the following:

\medskip
{\bf Theorem 4.2.} {\em $(1)$ Any theory $T_0$ of unary predicates
$P_i$, $i\in I$, and with infinite models, is exactly
$|T_0|$-approximable by the class $\mathcal{T}_0$ of theories of
unary predicates, in finite languages and with finite models.

$(2)$ Any theory $T_0$ of unary predicates $P_i$, $i\in I$, and
with infinite models, is countably approximable, by an appropriate
class $\mathcal{T}_0$ of theories of unary predicates.

$(3)$ A theory $T_1$ of unary predicates and with finite models is
{\rm (}countably{\rm )} approximable {\rm (}by an appropriate
class{\rm )} if and only if $T_1$ has an infinite language.}

\medskip
Proof. (1) Since the theory $T_0$ is based by the formulas
describing cardinality estimations for intersections of unary
predicates $P_i$, $i\in I$, and their complements, i.e.,
$P^{\delta_1}_{i_1}\wedge\ldots\wedge P^{\delta_k}_{i_k}$,
$\delta_j\in\{0,1\}$, then each formula $\varphi\in T_0$ belongs
to some theory in $\mathcal{T}_0$ whose models satisfy these
cardinality estimations. Since $|T_0|=|I|+\omega$, $T_0$ is
$|T_0|$-approximable by $\mathcal{T}_0$. Since the theories in
$\mathcal{T}_0$ have finite languages $T_0$ can not be
$\mu$-approximable for $\mu<|T_0|$, i.e., $T_0$ is exactly
$|T_0|$-approximable.

(2) If the language for $T_0$ is at most countable then we apply
the item 1. The same approach is valid if $T_0$ has at most
countably many independent predicates, i.e., all predicates are
boolean combinations of some at most countable family of them. So
below we assume that the language $\Sigma(T_0)$ is infinite and,
moreover, there are uncountably many independent predicates.

If $T_0$ has an infinite and co-infinite predicate $P_i$ then we
approximate links of $P_i$ with $\Sigma_i=\{P_j\mid j\in
I\setminus\{i\}\}$ by a countable family $\mathcal{T}_i$ of
theories $T'_k\ne T_0$, $k\in\omega$, in the language $P_i$, $i\in
I$, such that
$T'_k\upharpoonright\Sigma_i=T_0\upharpoonright\Sigma_i$ and $P_i$
for $T'_k$ corresponds to $T_0$ step-by-step, with respect to some
countable strictly increasing chain $(\Sigma'_k)_{k\in\omega}$,
where $\Sigma_i=\bigcup\limits_{k\in\omega}\Sigma'_k$.

If $T_0$ contains just (co-)finite predicates, we put, for $T'_k$,
predicates in $\Sigma'_k\cup\{P_i\}$ as required, and
$P\in\Sigma_i\setminus\Sigma'_k$ is either empty or complete such
that $P$ is empty for $T'_k$ if and only if $P$ is not empty for
$T_0$.

(3) If $T_1$ has a finite language then $T_1$ is isolated by a
sentence describing all connections between elements in a model of
$T_1$. Otherwise we can approximate $T_1$ by a countable family,
as in the item 2.~$\Box$

\medskip
{\bf Definition} \cite{Mak1}. A consistent sentence $\varphi$ of a
language $\Sigma$ is called {$\Sigma$-complete}, or simply {\em
complete}, if $\varphi$ forces a complete theory of the language
$\Sigma$.

\medskip
Clearly, a theory $T$ of a finite language $\Sigma$ contains a
$\Sigma$-complete sentence if and only if $T$ is finitely
axiomatizable (by that sentence).

\medskip
{\bf Theorem 4.3.} {\em For any theory $T$ the following
conditions are equivalent:

$(1)$ $T$ is $\lambda$-approximable for some $\lambda$;

$(2)$ $T$ is $\omega$-approximable;


$(3)$ the language $\Sigma(T)$ is finite and $T$ does not contain
a $\Sigma(T)$-complete sentence, or $\Sigma(T)$ is infinite.}

\medskip
Proof. $(2)\Rightarrow(1)$ is obvious.

$(1)\Rightarrow(3)$. If $\Sigma(T)$ is finite then the conclusion
follows by Proposition 3.1.

$(3)\Rightarrow(2)$. Let $\Sigma(T)$ is finite, and by assumption
$T$ is not finitely axiomatizable. Then by Proposition 3.1 the
theory $T$ is approximable, and being countable $T$ is
$\omega$-approximable.

If $\Sigma(T)$ is infinite we approximate $T$ by a countable
family $\{T_n\mid n\in\omega\}$ as shown in the proof of Theorem
4.2: again considering a strictly increasing family of languages
$\Sigma_n$, $n\in\omega$,
$\Sigma(T)=\bigcup\limits_{n\in\omega}\Sigma_n$, and
$\Sigma(T)$-theories $T_n$ such that
$T_n\upharpoonright\Sigma_n=T\upharpoonright\Sigma_n$ and
$P\in\Sigma(T)\setminus\Sigma_n$ is either empty or complete,
where $P$ is empty/complete for $T_n$ if and only if $P$ is not
empty/complete for $T$.~$\Box$

\medskip
Any approximating family $\{T_n\mid n\in\omega\}$ for $T$ in the
proof of Theorem 4.3 is called {\em trivial}, or {\em standard}.

Thus by Theorem 4.3 each theory, in an infinite language, has a
standard approximating family.

\section{Approximations by almost language uniform theories}

\medskip {\bf Definition} \cite{ccct}. A theory $T$ in a predicate language
$\Sigma$ is called {\em almost language
uniform},\index{Theory!almost language uniform} or a {\em {\rm
ALU}-theory}\index{{\rm ALU}-theory} if
for each arity $n$ with $n$-ary predicates for $\Sigma$ there is a
partition for all $n$-ary predicates, corresponding to the symbols
in $\Sigma$, with finitely many classes $K$ such that
any substitution preserving these classes preserves $T$, too.

\medskip
Below we consider approximations of theories by families of {\rm
ALU}-theories. Theories with these approximations are called {\em
{\rm ALU}-approximable}.

Since theories in $\overline{\mathcal{T}}_{\rm fin}$, being
theories of finite structures, are {\rm ALU}-theories
(\cite[Proposition 5.1]{ccct}) then theories $T$, which are
approximable by families
$\mathcal{T}\subset\overline{\mathcal{T}}_{\rm fin}$ are {\rm
ALU}-approximable. In particular, by Theorem 4.2 (1), theories of
unary predicates, with infinite models, are {\rm
ALU}-approximable.

By the definition each theory in a finite language is an {\rm
ALU}-theory. Since theories in a family for an approximation
satisfy the required approximable theory step-by-step, some
standard approximating family $\{T_n\mid n\in\omega\}$ for a
countable theory $T$, in an infinite language, can consist of {\rm
ALU}-theories such that for each $T_n$ only finitely many
predicates can differ from empty or complete ones. Thus we have
the following:

\medskip
{\bf Proposition 5.1.} {\em Any countable theory is an {\rm
ALU}-theory or it is {\rm ALU}-approximable.}

\section{Approximating subfamilies}

\medskip
{\bf Theorem 6.1.} {\em A family $\mathcal{T}$ of theories
contains an approximating subfamily if and only if $\mathcal{T}$
is infinite.}

\medskip
Proof. Since any approximating family is infinite then, having an
approximating subfamily, $\mathcal{T}$ is infinite.

Conversely, let $\mathcal{T}$ be infinite. Firstly, we assume that
the language $\Sigma=\Sigma(\mathcal{T})$ of $\mathcal{T}$ is at
most countable. We enumerate all $\Sigma$-sentences: $\varphi_n$,
$n\in\omega$, and construct an accumulation point for
$\mathcal{T}$ by induction. Since $\mathcal{T}_{\varphi_0}$ or
$\mathcal{T}_{\neg\varphi_0}$ is infinite  we can choose
$\psi_0=\varphi^{\delta}_0$ with infinite
$\mathcal{T}_{\varphi^\delta_0}$, $\delta\in\{0,1\}$. If $\psi_n$
is already defined, with infinite $\mathcal{T}_{\psi_n}$, then we
choose $\psi_{n+1}=\psi_n\wedge\varphi^\delta_{n+1}$, with
$\delta\in\{0,1\}$, such that $\mathcal{T}_{\psi_{n+1}}$ is
infinite. Finally, the set $\{\psi_n\mid n\in\omega\}$ forces a
complete theory $T$ being an accumulation point both for
$\mathcal{T}$ and for each $\mathcal{T}_{\psi_n}$. Thus,
$\mathcal{T}\setminus\{T\}$ is a required approximating family.

If $\Sigma$ is uncountable we find an accumulation point $T_0$ for
infinite $\mathcal{T}\upharpoonright\Sigma_0$, where $\Sigma_0$ is
a countable sublanguage of $\Sigma$. Now we extend $T_0$ till a
complete $\Sigma$-theory $T$ adding $\Sigma$-sentences $\chi$ such
that $\mathcal{T}_\chi$ are infinite. Again
$\mathcal{T}\setminus\{T\}$ is a required approximating
family.~$\Box$

\medskip
Using the construction for the proof of Theorem 6.1 we observe
that having infinite $\mathcal{T}_\varphi$ we obtain an
accumulation point $T$ for $\mathcal{T}_\varphi$ such that
$\varphi\in T$. So having infinite $\mathcal{T}_\varphi$ and
$\mathcal{T}_{\neg\varphi}$ we have at least two accumulation
points for $\mathcal{T}$. Therefore we obtain the following:

\medskip
{\bf Proposition 6.2.} {\em If a family $\mathcal{T}$ has infinite
subfamilies $\mathcal{T}_\varphi$ and $\mathcal{T}_{\neg\varphi}$
then $e$-${\rm Sp}(\mathcal{T})\geq 2$.}

\medskip
Similarly we have the following:

\medskip
{\bf Proposition 6.3.} {\em If a family $\mathcal{T}$ has infinite
subfamilies $\mathcal{T}_{\varphi_i}$ for pairwise inconsistent
formulas $\varphi_i$, $i\in I$, then $e$-${\rm
Sp}(\mathcal{T})\geq |I|$.}

\section{Single-valued approximations}

{\bf Definition.} An approximating family $\mathcal{T}$ is {\em
single-valued}, or {\em $e$-categorical}, if $e$-${\rm
Sp}(\mathcal{T})=1$.

\medskip
Clearly, if $\mathcal{T}$ is single-valued then $\mathcal{T}$ has
a single accumulation point, i.e., approximating some unique
theory $T$.

\medskip
If $T$ is the (unique) accumulation point for $\mathcal{T}$ then
the family $\mathcal{T}\cup\{T\}$ is also called {\em
single-valued}, or {\em $e$-categorical}.

\medskip
{\bf Proposition 7.1.} {\em Any $E$-closed family $\mathcal{T}$
with finite $e$-${\rm Sp}(\mathcal{T})>0$ is represented as a
disjoint union of $e$-categorical families
$\mathcal{T}_1,\ldots,\mathcal{T}_n$.}

\medskip
Proof. Let $e$-${\rm Sp}(\mathcal{T})=n$ and $T_1,\ldots,T_n$ be
accumulation points for $\mathcal{T}$ witnessing that equality.
Now we consider pairwise inconsistent formulas $\varphi_i\in T_i$
separating $T_i$ from $T_j$, $j\ne i$, i.e., with
$\neg\varphi_i\in T_j$. By Proposition 1.1 each family
$\mathcal{T}_i=\mathcal{T}_{\varphi_i}$ is infinite, with unique
accumulation point $T_i$, and thus $\mathcal{T}_i$ is
$e$-categorical. Besides, the families $\mathcal{T}_i$ are
disjoint by the choice of $\varphi_i$, and
$\mathcal{T}'=\mathcal{T}\setminus\left(\bigcup\limits_{i=1}^n\mathcal{T}_i\right)$
does not have accumulation points. Therefore
$\mathcal{T}'\cup\mathcal{T}_1$ is $e$-categorical, too. Thus,
$\mathcal{T}'\cup\mathcal{T}_1,\mathcal{T}_2,\ldots,\mathcal{T}_n$
is the required partition of $\mathcal{T}$ on $e$-categorical
families.~$\Box$

\medskip
{\bf Remark 7.2.} An arbitrary partition of a family $\mathcal{T}$
into disjoint e-categorical families $\mathcal{T}_i$, $i\in I$,
does not imply $e$-${\rm Sp}(\mathcal{T})=|I|$. Indeed, taking a
language $\{Q^k_n\mid k=1,2, n\in\omega\}$ of unary predicates we
can form a family $\{T^k_n\mid k=1,2, n\in\omega\}$ of theories
$T^k_n$ such that the predicates $Q^k_m$, $m\geq n$, are complete
and the others are empty. All families $\mathcal{T}$, $\{T^1_n\mid
n\in\omega\}$, $\{T^2_n\mid n\in\omega\}$ are $e$-categorical,
whose common accumulation point consists of all empty predicates,
whereas $\{T^1_n\mid n\in\omega\}$, $\{T^2_n\mid n\in\omega\}$
form a partition of $\mathcal{T}$.

Similarly, having an arbitrary infinite family $\mathcal{T}$ of
theories in the empty language (which is $e$-categorical) we can
arbitrarily divide $\mathcal{T}$ into two infinite parts each of
which is again $e$-categorical, with the common accumulation point
$T$ having infinite models.

More generally, by Theorem 6.1, if $\mathcal{T}$ is
$e$-categorical then each infinite
$\mathcal{T}'\subseteq\mathcal{T}$ is again $e$-categorical with
the same accumulation point.

\medskip
{\bf Definition.} An approximating family $\mathcal{T}$ is called
{\em $e$-minimal} if for any sentence $\varphi\in\Sigma(T)$,
$\mathcal{T}_\varphi$ is finite or $\mathcal{T}_{\neg\varphi}$ is
finite.

\medskip
{\bf Theorem 7.3.} {\em A family $\mathcal{T}$ is $e$-minimal if
and only if it is $e$-categorical.}

\medskip
Proof. Let $\mathcal{T}$ be $e$-minimal. We consider the set
$T=\{\varphi\in\Sigma(\mathcal{T})\mid\mathcal{T}_\varphi\mbox{ is
infinite}\}$. By compactness $T$ is consistent and by
$e$-minimality of $\mathcal{T}$, $T$ is a complete theory. Thus,
by the definition $T$ is an accumulation point for $\mathcal{T}$.
This accumulation point is unique since if $T'\ne T$ is a
$\Sigma(\mathcal{T})$-theory then there is $\varphi\in T$ such
that $\neg\varphi\in T'$, so $\mathcal{T}_{\neg\varphi}$ and
$T'\notin{\rm Cl}_E(\mathcal{T}\setminus\{T'\})$ by Proposition
1.1. Thus, $\mathcal{T}$ is $e$-categorical.

Conversely, if $\mathcal{T}$ is not $e$-minimal then $e$-${\rm
Sp}(\mathcal{T})\geq 2$ by Corollary 6.2. Thus, $\mathcal{T}$ is
not $e$-categorical.~$\Box$

\medskip
{\bf Remark 7.4.} As shown in Remark 7.2 $e$-categorical families
can be always divided into $e$-categorical parts. So the condition
for $e$-minimality, on divisibilities only with respect to
neighbourhoods $\mathcal{T}_\varphi$, is essential counting
$e$-spectra.

\noindent Sobolev Institute of Mathematics, \\ 4, Acad. Koptyug
avenue, Novosibirsk, 630090, Russia; \\ Novosibirsk State
Technical
University, \\ 20, K.Marx avenue, Novosibirsk, 630073, Russia; \\
Novosibirsk State University, \\ 1, Pirogova street, Novosibirsk,
630090, Russia
\end{document}